\documentclass[a4paper]{article}
\usepackage[utf8]{inputenc}
\usepackage[T1,T2A,OT2,OT1]{fontenc}
\usepackage{amssymb}
\usepackage{amsmath}
\usepackage{mathrsfs}
\usepackage{amsfonts}
\usepackage{color}
\usepackage{vmargin}
\usepackage{amsthm}
\usepackage{graphicx}
\usepackage{xfrac}
\usepackage[sort&compress, numbers]{natbib}
\usepackage[russian,english]{babel}

\setmarginsrb{20mm}{20mm}{20mm}{20mm}{10mm}{10mm}{10mm}{10mm}

\newtheoremstyle{remark}
  {}{}{}{}{\bfseries}{.}{.5em}{{\thmname{#1 }}{\thmnumber{#2}}{\thmnote{ (#3)}}}

\RequirePackage{amsthm}

\newtheorem{defi}{Definition}[section]

\newtheorem{tw}[defi]{Theorem}

\newtheorem{cor}[defi]{Corollary}
\newtheorem{prop}[defi]{Proposition}

\def\={\hspace{-3mm}&=&\hspace{-3mm}}

\renewcommand{\r}{\mathbb{R}}
\newcommand{\n}{\mathbb{N}}

\newcommand{\ioi}{\int_0^{\infty}}
\newcommand{\io}{\int_{\Omega}}
\newcommand{\hoo}{H^1_0(\Omega)}
\newcommand{\ld}{L^2(\Omega)}
\newcommand{\loc}{\textnormal{loc}}
\newcommand{\ul}{\frac{1}{2}}

\let\phi\varphi
\let\epsilon\varepsilon

\begin{document}

\title{\bf Decay of solutions of non-homogenous hyperbolic equations}

\author{Piotr Micha{\l} Bies\\
\it\small{Faculty of Mathematics and Information Sciences,}\\
\it\small{Warsaw University of Technology,}\\
\it\small{Ul. Koszykowa 75, 00-662 Warsaw, Poland.}\\
{\tt biesp@mini.pw.edu.pl}}
\maketitle

\begin{abstract}
We consider conditions for the decay in time of solutions of non-homogenous hyperbolic equations. It is proven that solutions of the equations go to $0$ in $L^2$ at infinity if and only if an equation's right-hand side uniquely determines the initial conditions in a certain way. We also obtain that a hyperbolic equation has a unique solution that fades when $t\to\infty$.
\end{abstract}
\maketitle
\bigskip

\medskip


\medskip
\section{Introduction}

Let $\Omega \subset\r^n$ be an open and bounded set. In the paper, we consider the decay in time of solutions to the following problem
\begin{equation}\label{wav}
\begin{cases}
u_t+Lu=f&\textrm{in }(0,\infty)\times\Omega,\\
u=0&\textrm{on }(0,\infty)\times\partial\Omega,\\
u=g, u_t=h&\textrm{on }\{0\}\times\Omega,
\end{cases}
\end{equation}
where functions $f,g,h$ are given and we want to find is $u$. We also require that operator $u_t+Lu$ is a second-order hyperbolic operator. In our considerations, the operator $L$ has the following form
\begin{align*}
Lu=-\sum_{i,j=1}^n\left(a^{ij}u_{x_i}\right)_{x_j},
\end{align*}
where we also assume that $a^{ij}$ depends only on $x\in\Omega$. Moreover, we require that there exists $\lambda>0$ such that 
\begin{align}\label{lamb}
\sum^n_{i,j=1}a^{ij}(x)y_iy_j\geq\lambda|y|^2 \textrm{ for all }x\in\Omega\textrm{ and }y\in\r^n.
\end{align}
In addition, we assume that $a^{ij}\in L^{\infty}(\Omega)$ and $a^{ij}=a^{ji}$ for $i,j=1,\ldots,n$.  

It is well known that $u$ can be interpreted as a displacement of a material. The coefficients $a^{ij}$ come from material principles. In general, equation \eqref{wav} models wave transmission in the medium. Problems as in \eqref{wav} were widely studied. There are well-known results about the existence, uniqueness, and regularity of solutions of hyperbolic equations (see \cite{Evans, Lady, Lions, Wloka}). 

The following theorem is the main result of the present article.
\begin{tw}\label{maintw}
Let us assume that $g\in\hoo,h\in\ld,f\in L^2_{\loc}\left([0,\infty);\ld\right)\cap L^1(0,\infty;\ld)$ and $u$ is a weak solution of the problem \eqref{wav}. Then, 
\begin{align}\label{tez1}
u(t,\cdot)\to 0 \textrm{ in }\ld\textrm{ when } t\to\infty
\end{align}
if and only if the following equalities
\begin{align}\label{tez2}
\begin{split}
(g,\phi_m)_{\ld}&=\frac1{\sqrt{\lambda_m}}\int_0^{\infty}\sin\left(\sqrt{\lambda_m} s\right) (f(s,\cdot),\phi_m)_{\ld}ds,\\ 
 (h,\phi_m)_{\ld}&=-\int_0^{\infty}\cos\left({\sqrt{\lambda_m}} s\right) (f(s,\cdot),\phi_m)_{\ld}ds
 \end{split}
\end{align}
hold for all $m\in\n$.
\end{tw}
We see that the conditions in \eqref{tez2} uniquely determine the initial conditions of the problem \eqref{wav} if we know that \eqref{tez1} is satisfied. The equalities uniquely connect the initial conditions with $f$ in \eqref{wav}.

The decay of wave phenomena is natural and often observed. We see it in the water, string, and many other cases. It also occurs in the hyperbolic-parabolic system of thermoelasticity. It is a problem that describes oscillations and the heat in a medium. It can be written as follows
\begin{equation}\label{eq}
\begin{cases}
u_{tt}-\Delta u=\mu\operatorname{div}\left(\theta I\right),&\textrm{in }(0,\infty)\times \Omega,\\
\theta_t-\Delta\theta=\mu\theta \operatorname{div}(u_t),&\textrm{in }(0,\infty)\times \Omega,\\
u|_{\partial\Omega}=0,\ \frac{\partial\theta}{\partial n}|_{\partial\Omega}=0&\textrm{on } (0,\infty),\\
u(0,\cdot)=u_0,\ u_t(0,\cdot)=v_0,\ \theta(0,\cdot)=\theta_0>0,
\end{cases}
\end{equation}
where $\mu$ is a constant, initial data $\theta_0, u_0, v_0$ are given and we want to find $u\colon[0,\infty)\times\Omega\to\r^n$, which is the displacement, and $\theta\colon[0,\infty)\times\Omega\to\r$, which is the temperature. It is usually considered for $n=1,2,3$. We see that the equation for $u$ is the non-homogenous wave equation. System \eqref{eq} and similar were widely studied in many cases (see \cite{BC, TC1, TC2, Racke1, Racke2, Racke3, MR0629700, MR1031681, MR1057652, MR4432953, MR0860899, MR1270660} and many others).

In our considerations, it is important that, by the second law of thermodynamics, we predict that 
\begin{align*}
\lim_{t\to\infty}u(t,\cdot)=0\quad \textrm{ and } \quad\lim_{t\to\infty}\theta(t,\cdot)=\textnormal{constant}.
\end{align*}
It has been shown lately in paper \cite{BCst} for $n=1$. It means that Theorem \ref{maintw} perhaps can be applied to the upper equation in \eqref{eq}.

The paper is divided into three sections. We consider the non-homogenous harmonic oscillator equation in Section \ref{sechar}. It is a well-known ordinary differential equation of the second order. It can be said that this is an ordinary version of a hyperbolic equation. It is shown that the solution of the harmonic oscillator equation disappears when $t\to\infty$ if and only if the initial values and the right-hand side satisfy certain conditions. They uniquely connect them. As a corollary, we obtain that the non-homogenous oscillator equation with a fading solution at infinity has a unique solution. 
In the last section, we consider a non-homogenous hyperbolic equation. There, we prove the main theorem. Again, it occurs that the equation has a unique solution in the class of functions that disappear when time goes to infinity.

\section{Harmonic oscillator}\label{sechar}

In this section, we consider an ordinary differential equation
\begin{equation}\label{osceq}
y''(t)+\mu^2 y(t)=f(t)\quad\textrm{ for }t\geq 0,
\end{equation}
where the constant $\mu>0$ and the function $f\colon [0,\infty)\to\r$ are given.
\begin{tw}\label{hartw}
Let $f\colon [0,\infty)\to\r$ be such that the integrals
\begin{equation}\label{zal1}
\ioi\sin(\mu s)f(s)ds,\quad \ioi\cos(\mu s)f(s)ds
\end{equation} 
are convergent in the Riemann sense and let $y$ satisfy equation \eqref{osceq}. Then,
$$
\lim_{t\to \infty}y(t)=0
$$
if and only if 
\begin{align}\label{inithar}
y(0)=\frac1\mu\ioi\sin(\mu s)f(s)ds,\quad y'(0)=-\ioi\cos(\mu s) f(s)ds.
\end{align}
\end{tw}
\begin{proof}
We know that the general solution of equation \eqref{osceq} is given by the formula
\begin{align*}
y(t)=\cos(\mu t)\left(C_1-\frac1\mu\int_0^t\sin(\mu s)f(s)ds\right)+\sin(\mu t)\left(C_2+\frac1\mu\int_0^t\cos(\mu s)f(s)ds\right)
\end{align*}
for $t\geq 0$ and for arbitrary $C_1,C_2\in\r$. 

Thus, we see that the solution of \eqref{osceq} with initial values \eqref{inithar} is as follows
\begin{align}\label{harfor}
\begin{split}
y(t)&=\cos(\mu t)\left(\frac1\mu\ioi\sin(\mu s) f(s)ds-\frac1\mu\int_0^t\sin(\mu s) f(s)ds\right)
\\&\quad +\sin(\mu t)\left(-\frac1\mu\ioi\cos(\mu s) f(s)ds+\frac1\mu\int_0^t\cos(\mu s)f(s)ds\right).
\end{split}
\end{align}
Hence, we have that $y(t)\to 0$, when $t\to \infty$. Therefore, if equalities \eqref{inithar} are satisfied, then $y(t)\to 0$, when $t\to \infty$.

Let us assume that one of the equalities from \eqref{inithar} is not satisfied. Then, we have that
\begin{align*}
\left(C_1-\frac1\mu\int_0^t\sin\mu sf(s)ds\right)\not\to 0\textrm{ or }\left(C_2+\frac1\mu\int_0^t\cos\mu sf(s)ds\right)\not\to 0, \textrm{ when }t\to\infty.
\end{align*}
Therefore, we see that the solution cannot $y(t)\to 0$ if $t\to\infty$.
\end{proof}

We finish the section with the following corollary.
\begin{cor}\label{corhar}
Let us assume that $f$ is such that the integrals in \eqref{zal1} are convergent in the Riemann sense, then the problem
\begin{align*}
\begin{cases}
y''(t)+\mu^2 y(t)=f(t)&\textrm{for }t\geq 0,\\
\lim\limits_{t\to\infty}y(t)=0
\end{cases}
\end{align*}
has a unique solution.
\end{cor}

\section{Wave equation}

In this section, we consider a second-order non-homogenous hyperbolic equation. We prove the main theorem of the paper here. Let $\{\phi_m\}$ be a sequence of the operator $L$ eigenfunctions in $\hoo$. We assume that $\|\phi_m\|_{\ld}=1$. We know (see for instance \cite{Evans}) that $\operatorname{span}\{\phi_m\}$ is a dense set in $\hoo$. Let $\{\lambda_m\}$ be a sequence of the operator $L$ eigenvalues. We suppose also that $g\in\hoo,h\in\ld$ and $f\in L^2_{\loc}\left([0,\infty);\ld\right)$. Let us remind here that with $L$ is associated the bilinear form
\begin{align*}
B(u,v)=\io\sum_{i,j=1}^na^{ij}u_{x_i}u_{x_j}dx.
\end{align*}

By the separation of variables technique, we know that the solution of the problem \eqref{wav} is given by the formula 
\begin{align}\label{szer1}
u(x,t)=\sum_{m=1}^{\infty}d_m(t)\phi_m(x)\textrm{ in }L^{\infty}_{\loc}([0,\infty),\ld)\cap C([0,\infty),\ld),
\end{align}
where functions $\{d_m\}$ satisfy
\begin{align}\label{equdm}
\begin{cases}
d_m''+\lambda_m d_m=(f,\phi_m)_{\ld}&\textrm{on }[0,\infty),\\
d_m(0)=(g,\phi_m)_{\ld},\ d_m'(0)=(h,\phi_m)_{\ld}
\end{cases}
\end{align}
for all $m\in\n$. We shall prove that the function in \eqref{szer1} is well defined.
\begin{prop}
Let us assume that $g\in\hoo,h\in\ld$ and $f\in L^2_{\loc}\left([0,\infty);\ld\right)$. Then, the function in formula \eqref{szer1} is well defined.
\end{prop}
\begin{proof}
Let us take $T>0$. We see that
\begin{align*}
u_m:=\sum_{k=1}^md_k\phi_k
\end{align*}
is a sequence from the Galerkin method applied to the problem \eqref{wav}. Thus, by \cite{Evans}, we know that $\{u_m\}$ is bounded in $L^{\infty}\left(0,T;\hoo\right)\cap W^{1,\infty}\left(0,T;\ld\right)\cap H^2(0,T;H^{-1}\left(\Omega)\right)$. Hence, we can take a subsequence of $\{u_m\}$ (still denoted as $\{u_m\}$) and $u\in L^{\infty}\left(0,T;\hoo\right)\cap W^{1,\infty}\left(0,T;\ld)\right)\cap H^2(0,T;H^{-1}\left(\Omega)\right)$ such that
\begin{align}\label{sl}
\begin{split}
u_m&\stackrel{*}{\rightharpoonup} u \textrm{ in }L^{\infty}\left(0,T;\hoo\right),\\
u_m'&\stackrel{*}{\rightharpoonup} u' \textrm{ in }L^{\infty}\left(0,T;\ld\right),\\
u_m''&\,{\rightharpoonup}\; u'' \textrm{ in }L^2\left(0,T;H^{-1}(\Omega)\right).
\end{split}
\end{align}
We know that $u$ is a solution of \eqref{wav}, which is unique. Therefore, in \eqref{sl}, we can take a whole sequence, not only the subsequence. By Aubin-Lions lemma (see \cite{Boyer}), we get that 
\begin{align*}
u_m\to u\textrm{ in }C\left([0,T],\ld\right).
\end{align*}
\end{proof}

Because $d_m$ satisfies \eqref{equdm}, so we have an analogous formula as in \eqref{harfor}
\begin{align}\label{ford}
\begin{split}
d_m(t)&=\cos\left(\sqrt{\lambda_m} t\right)\left((g,\phi_m)_{\ld}-\frac1{\sqrt{\lambda_m}}\int_0^t\sin\left(\sqrt{\lambda_m} s\right) (f(s,\cdot),\phi_m)_{\ld}ds\right)
\\&\quad +\sin\left(\sqrt{\lambda_m} t\right)\left(\frac1{\sqrt{\lambda_m}}(h,\phi_m)_{\ld}+\frac1{\sqrt{\lambda_m}}\int_0^t\cos\left({\sqrt{\lambda_m}} s\right) (f(s,\cdot),\phi_m)_{\ld}ds\right)
\end{split}
\end{align}
for all $m\in\n$.

Now, we can prove the paper's main result.
\begin{proof}[Proof of Theorem \ref{maintw}]
Let us assume that \eqref{tez1} is satisfied, and let us take $m\in\n$. Then, $(u(t,\cdot),\phi_m)_{\ld}=d_m(t)\to 0$, when $t\to\infty$. By Theorem \ref{hartw}, it implies \eqref{tez2}.

Next, let us assume that \eqref{tez2} is satisfied for all $m\in\n$. We need a certain inequality before we will prove the thesis. First, we will justify it for $$u_m:=\sum^m_{k=1}d_k\phi_k,$$
where $d_k$ are defined above.

Let us take $u_m'$ as a test function in the weak formulation of the problem for $u_m$. We get
\begin{align*}
\io u_{m}''u_m'dx+B(Du_m, Du_m')=\io fu_m'dx.
\end{align*}
It simply yields
\begin{align*}
\frac12\frac d{dt}\left(\io u_{m,t}^2dx+B(Du_m,Du_m)\right)=\io fu_m'dx.
\end{align*}
Now, Let us fix $T>0$ and let us take $t\in [0,T]$. We integrate the above equality over $[0,t]$. It gives us
\begin{align*}
\io u_{m,t}^2dx+B(u_m,u_m)&=2\int_0^t\io fu_t\,dx\, ds+\|h\|^2_{\ld}+B(g,g)\\
&\leq 2\|u_m'\|_{L^{\infty}(0,T;\ld)}\|f\|_{L^1(0,T;\ld)}+\|h\|^2_{\ld}+C\|g\|^2_{\hoo}\\
&\leq \frac14\|u_m'\|_{L^{\infty}(0,T;\ld)}^2+4\|f\|_{L^1(0,\infty;\ld)}^2+\|h\|^2_{\ld}+C\|g\|^2_{\hoo},
\end{align*}
where $C=(\sum_{i,j=1}^n\|a^{ij}\|_{L^{\infty}(\Omega)})^{\ul}$.
It implies that 
\begin{align*}
\frac12&\left(\|u'_m\|_{L^{\infty}(0,T;\ld)}^2+\lambda\|u_m\|_{L^{\infty}(0,T;\hoo)}^2\right)\leq\sup_{t\in [0,T]}\left(\io u_{m,t}^2dx+\lambda\io|Du_m|^2dx\right)\\
&\leq \frac14\|u_m'\|_{L^{\infty}(0,T;\ld)}^2+4\|f\|_{L^1(0,\infty;\ld)}^2+\|h\|^2_{\ld}+C\|g\|^2_{\hoo},
\end{align*}
where $\lambda$ is defined in \eqref{lamb}.
Hence, we arrive with
\begin{align*}
\frac\lambda2\|u_m\|_{L^{\infty}(0,T;\hoo)}^2&\leq \frac14\|u_m'\|_{L^{\infty}(0,T;\ld)}^2+\frac\lambda2\|u_m\|_{L^{\infty}(0,T;\hoo)}^2\\
&\leq4\|f\|_{L^1(0,\infty;\ld)}^2+\|h\|^2_{\ld}+\|g\|^2_{\hoo}.
\end{align*}
Because $T>0$ was arbitrary, so we see that the sequence is bounded $\{u_m\}$ in $L^{\infty}(0,\infty;\hoo)$. Thus, because $u$ is a *-weak limit of $\{u_m\}$, it derives that $u\in L^{\infty}(0,\infty;\hoo)$.

Now, we think about $u$ as representative such that $u\in C\left([0,\infty),\ld\right)$. Moreover, let us take a measurable set $X\subset[0,\infty)$ such that $\left|[0,\infty)\setminus X\right|=0$ and the function $u|_X$ is bounded as a function with values in $\hoo$. Next, let us take a sequence $\{t_k\}\subset X$ such that $t_k\to\infty$. A sequence $\{u(t_k,\cdot)\}$ is bounded in $\hoo$. Moreover for an arbitrary $m\in\n$, we have 
\begin{align*}
(u(t_k,\cdot),\phi_m)_{\hoo}=d_m(t_k)\|\phi_m\|_{\hoo}^2=d_m(t_k)\lambda_m\to 0\textrm{ when } k\to\infty.
\end{align*}
Because $\operatorname{span}\{\phi_m\}$ is dense in $\hoo$, so we obtain 
\begin{align*}
u(t_k,\cdot)\rightharpoonup 0\textrm{ in }\hoo.
\end{align*}
By the compact embedding, we have that
\begin{align*}
u(t_k,\cdot)\to 0\textrm{ in }\ld.
\end{align*}
Therefore, we obtain
\begin{align*}
 u(t,\cdot)\to 0 \textrm{ in }\ld\textrm{ when }t\in X \textrm{ and }t\to \infty.
\end{align*}

Because $u\in C\left([0,\infty),\ld\right)$, so one can easy prove that 
\begin{align*}
u(t,\cdot)\to 0 \textrm{ in } \ld \textrm{ when }t\to\infty.
\end{align*}
Details are left to the reader.
\end{proof}

Ultimately, we formulate and prove an analogous result as Corollary \ref{corhar} for the wave equation.
\begin{prop}
Let us assume that the hypotheses of Theorem \ref{maintw} are satisfied. Then, the problem
\begin{align}\label{decaypr}
\begin{cases}
u_{tt}+Lu=f&\textrm{in }(0,\infty)\times\Omega,\\
u=0& \textrm{on } (0,\infty)\times\partial\Omega,\\
\lim\limits_{t\to\infty}\|u(t,\cdot)\|_{\ld}=0,
\end{cases}
\end{align}
has a unique solution in $L^{\infty}_{\loc}\left(0,\infty;\hoo\right)\cap W^{1,\infty}_{\loc}\left(0,\infty;\ld\right)\cap H^2_{\loc}(0,\infty;H^{-1}\left(\Omega)\right)$.
\end{prop}
\begin{proof}
We must show that quantities in \eqref{tez2} define unique initial values for $u$. Let us remind that $\{\phi_m\}$ are the eigenfunctions of $L$ and $\{\lambda_m\}$ are eigenvalues of it. The functions $\{\phi_m\}$ are also normalized in $\ld$. Let us denote
\begin{align*}
f_m(t):=(f(t,\cdot),\phi)_{\ld}
\end{align*}
for all $t\in[0,\infty)$ and $m\in\n$. Then, we have
\begin{align*}
f(t,x):=\sum_{m=1}^{\infty}f_m(t)\phi_m(x)
\end{align*}
for almost all $(t,x)\in [0,\infty)\times\Omega$ and series is convergent in $\ld$ and pointwise for almost all $t\in[0,\infty)$. It also derives
\begin{align*}
\ioi\|f(t,\cdot)\|_{\ld}dt=\ioi\left(\sum_{m=1}^{\infty}f_m^2(t)\right)^{\ul}dt
\end{align*}
and the integrals are finite because $f\in L^1(0,\infty;\ld)$.

First, we will show that the formula in \eqref{tez2} can well define the velocity. Let us denote
\begin{align*}
h_m=-\int_0^{\infty}\cos\left({\sqrt{\lambda_m}} s\right) (f(s,\cdot),\phi_m)_{\ld}ds
\end{align*}
for all $m\in\n$. We want to show that $\{h_m\}\in\ell^2$. The Minkowski generalized inequality (see \cite{Fourier})implies
\begin{align*}
\sum_{m=1}^{\infty}h_m^2=\sum_{m=1}^{\infty}\left(\int_0^{\infty}\cos\left({\sqrt{\lambda_m}} s\right) (f(s,\cdot),\phi_m)_{\ld}ds\right)^2
\leq\left(\ioi\left(\sum_{m=1}^{\infty}f_m^2(t)\right)^{\ul}dt\right)^2.
\end{align*}
 The integral on the right-hand side is finite, so $\{h_m\}\in\ell^2$. Therefore, the formula
 \begin{align*}
 h:=\sum_{m=1}^{\infty}h_m\phi_m
 \end{align*}
 defines the function $h$ properly.
 
Next, we are going to the initial value for $u$. Henceforth, we will write
\begin{align*}
\|w\|_{\hoo}=\left(B(w,w)\right)^{\ul}
\end{align*}
for $w\in \hoo$. It is equivalent norm to the standard norm in $\hoo$. Let us denote 
\begin{align*}
g_m:=\frac1{\sqrt{\lambda_m}}\int_0^{\infty}\sin\left(\sqrt{\lambda_m} s\right) (f(s,\cdot),\phi_m)_{\ld}ds.
\end{align*}
One can easily see that
\begin{align*}
\|\phi_m\|_{\hoo}=\sqrt{\lambda_m}.
\end{align*}
On the other hand, we show that $\{h_m\sqrt{\lambda_m}\}\in\ell^2$ similarly as above. Therefore, the formula
\begin{align*}
g:=\sum_{m=1}^{\infty}g_m{\phi_m}
\end{align*}
defines $g\in \hoo$. We see that $f,g,h$ satisfy condition \eqref{tez2}.

Let us consider $u\in L^{\infty}_{\loc}\left(0,\infty;\hoo\right)\cap W^{1,\infty}_{\loc}\left(0,\infty;\ld\right)\cap H^2_{\loc}(0,\infty;H^{-1}\left(\Omega)\right)$ the weak solution of the problem
\begin{equation}\label{pompr}
\begin{cases}
u_t+Lu=f&\textrm{in }(0,\infty)\times\Omega,\\
u=0&\textrm{on }(0,\infty)\times\partial\Omega,\\
u=g, u_t=h&\textrm{on }\{0\}\times\Omega.
\end{cases}
\end{equation}
Then, Theorem \ref{maintw} gives us that 
\begin{align*}
\lim_{t\to\infty}\|u(t,\cdot)\|_{\ld}=0.
\end{align*}
On the other, let $w$ be another weak solution of \eqref{decaypr}. Then, by Theorem \ref{maintw}, it has to be a weak solution of \eqref{pompr}. Because it is unique, so $w=u$.
\end{proof}
\subsection*{Acknowledgments}
P. M. Bies would like to thank Tomasz Cie\'slak for reading the first version of the manuscript and for general help during the paper preparation. The author also wishes to thank Professor Krzysztof Che\l mi\'nski and Professor Przemys\l aw G\'orka for finding gaps in the proof of Theorem \ref{maintw}.
\bibliographystyle{abbrv}
\bibliography{WaveRight}

\begin{thebibliography}{10}

\bibitem{Fourier}
H.~Bahouri, J.-Y. Chemin, and R.~Danchin.
\newblock {\em Fourier analysis and nonlinear partial differential equations},
  volume 343 of {\em Grundlehren der mathematischen Wissenschaften [Fundamental
  Principles of Mathematical Sciences]}.
\newblock Springer, Heidelberg, 2011.

\bibitem{BC}
P.~M. Bies and T.~Cie\'{s}lak.
\newblock Global-in-time regular unique solutions with positive temperature to
  one-dimensional thermoelasticity.
\newblock {\em SIAM J. Math. Anal.}, 55(6):7024--7038, 2023.

\bibitem{BCst}
P.~M. Bies and T.~Cieślak.
\newblock Time-asymptotics of a heated string.
\newblock 2024.
\newblock arXiv:2405.04310.

\bibitem{Boyer}
F.~Boyer and P.~Fabrie.
\newblock {\em Mathematical tools for the study of the incompressible
  {N}avier-{S}tokes equations and related models}, volume 183 of {\em Applied
  Mathematical Sciences}.
\newblock Springer, New York, 2013.

\bibitem{MR1270660}
G.~Q. Chen and C.~M. Dafermos.
\newblock The vanishing viscosity method in one-dimensional thermoelasticity.
\newblock {\em Trans. Amer. Math. Soc.}, 347(2):531--541, 1995.

\bibitem{TC2}
T.~Cie\'{s}lak, M.~Gali\'{c}, and B.~Muha.
\newblock A model in one-dimensional thermoelasticity.
\newblock {\em Nonlinear Anal.}, 216:Paper No. 112703, 21, 2022.

\bibitem{TC1}
T.~Cie\'{s}lak, B.~Muha, and S.~a. Trifunovi\'{c}.
\newblock Global weak solutions in nonlinear 3{D} thermoelasticity.
\newblock {\em Calc. Var. Partial Differential Equations}, 63(1):Paper No. 26,
  36, 2024.

\bibitem{MR0860899}
C.~M. Dafermos and L.~Hsiao.
\newblock Development of singularities in solutions of the equations of
  nonlinear thermoelasticity.
\newblock {\em Quart. Appl. Math.}, 44(3):463--474, 1986.

\bibitem{Evans}
L.~C. Evans.
\newblock {\em Partial differential equations}, volume~19 of {\em Graduate
  Studies in Mathematics}.
\newblock American Mathematical Society, Providence, RI, second edition, 2010.

\bibitem{MR4432953}
M.~Galanopoulou, A.~Vikelis, and K.~Koumatos.
\newblock Weak-strong uniqueness for measure-valued solutions to the equations
  of quasiconvex adiabatic thermoelasticity.
\newblock {\em Comm. Partial Differential Equations}, 47(6):1133--1175, 2022.

\bibitem{MR1057652}
W.~J. Hrusa and S.~A. Messaoudi.
\newblock On formation of singularities in one-dimensional nonlinear
  thermoelasticity.
\newblock {\em Arch. Rational Mech. Anal.}, 111(2):135--151, 1990.

\bibitem{MR1031681}
W.~J. Hrusa and M.~A. Tarabek.
\newblock On smooth solutions of the {C}auchy problem in one-dimensional
  nonlinear thermoelasticity.
\newblock {\em Quart. Appl. Math.}, 47(4):631--644, 1989.

\bibitem{Lady}
O.~A. Ladyzhenskaya.
\newblock {\em {\cyr Kraevye zadachi matematichesko\u{\i} fiziki}}.
\newblock Izdat. ``Nauka'', Moscow, 1973.

\bibitem{Lions}
J.-L. Lions and E.~Magenes.
\newblock {\em Non-homogeneous boundary value problems and applications. {V}ol.
  {II}}, volume Band 182 of {\em Die Grundlehren der mathematischen
  Wissenschaften}.
\newblock Springer-Verlag, New York-Heidelberg, 1972.
\newblock Translated from the French by P. Kenneth.

\bibitem{Racke1}
R.~Racke.
\newblock Initial boundary value problems in one-dimensional nonlinear
  thermoelasticity.
\newblock {\em Math. Methods Appl. Sci.}, 10(5):517--529, 1988.

\bibitem{Racke2}
R.~Racke and Y.~Shibata.
\newblock Global smooth solutions and asymptotic stability in one-dimensional
  nonlinear thermoelasticity.
\newblock {\em Arch. Rational Mech. Anal.}, 116(1):1--34, 1991.

\bibitem{Racke3}
R.~Racke, Y.~Shibata, and S.~M. Zheng.
\newblock Global solvability and exponential stability in one-dimensional
  nonlinear thermoelasticity.
\newblock {\em Quart. Appl. Math.}, 51(4):751--763, 1993.

\bibitem{MR0629700}
M.~Slemrod.
\newblock Global existence, uniqueness, and asymptotic stability of classical
  smooth solutions in one-dimensional nonlinear thermoelasticity.
\newblock {\em Arch. Rational Mech. Anal.}, 76(2):97--133, 1981.

\bibitem{Wloka}
J.~Wloka.
\newblock {\em Partial differential equations}.
\newblock Cambridge University Press, Cambridge, 1987.
\newblock Translated from the German by C. B. Thomas and M. J. Thomas.

\end{thebibliography}
\end{document}